\documentclass{article}%
\usepackage{amsmath}%
\usepackage{amsfonts}%
\usepackage{amssymb}%
\usepackage{graphicx}

\begin{document}

\title{On a Method for Solving Infinite Series}
\author{Henrik Stenlund\thanks{The author is obliged to Visilab Signal Technologies for supporting this work.}
\\Visilab Signal Technologies Oy, Finland}
\date{6th May, 2014}
\maketitle
\begin{abstract}
This paper is about a method for solving infinite series in closed form by using inverse and forward Laplace transforms. The resulting integral is to be solved instead. The method is extended by parametrizing the series. A further Laplace transform with respect to it will offer more options for solving a series.
 \footnote{Visilab Report \#2014-05} \footnote{v2 two references added with comments}

\subsection{Keywords}
Summation of series, infinite series, Laplace transform, inverse Laplace transform
\subsection{Mathematical Classification}
MSC: 44A10, 11M41, 16W60, 20F14, 40A25, 65B10
\end{abstract}

There is a crack in the mountain.

\tableofcontents

\section{Introduction}
\subsection{General}
An infinite series of the form below, assuming it to converge,
\begin{equation}
a=\sum_{k=1}^{\infty}{g(k)} \label{eqn1}
\end{equation}
needs to be solved. If there are any parameters involved, the results are either in a closed form or in terms of some special function which may itself be defined as an infinite series. Else the result is purely numerical but may consist of a function of constants, like $\pi$. Another infinite series or an integral is often acceptable as a solution being more suitable for further analysis. There exists just a handful of general ways for series solutions. Of those the most important ones are Euler method \cite{Euler1755}, the Euler-Maclaurin and the Poisson formulas \cite{Abramowitz1970}, \cite{Ivic1985}, \cite{Patterson1995}, in addition to using integral transforms. Other very specific methods exist \cite{Whittaker1915} but they have less general interest. Such are the Ramanujan summation \cite{Candelpergher2009} and Voronoi summation \cite{Ivic1985}. Hardy's monograph \cite{Hardy1991} on diverging series is most useful showing ways of handling various types of series. 

Wheelon \cite{Wheelon1954} has developed a method for solving series in closed form. He starts from a known Laplace transform pair and creates a sum to get an integral. The method resembles the simple one described in this paper though the approach is different. MacFarlane \cite{MacFarlane1949} has an approach analogous to the previous by using the Mellin transform. McFadden \cite{McFadden1953} and Glasser \cite{Glasser1971} have shared interesting formulas for solving series. The latter ones are similar to the simple method here. 

\subsection{Euler Method}
Euler \cite{Euler1755} has given a method which can sometimes be applied to solve simple convergent series of type
\begin{equation}
e=\sum_{k=1}^{\infty}{a_k} \label{eqn4}
\end{equation}
The series in equation (\ref{eqn4}) can be extended to a power series as follows
\begin{equation}
e(x)=\sum_{k=1}^{\infty}{a_k{x^k}} \label{eqn6}
\end{equation}
This power series can possibly be identified as some known expansion. Then the limit $x\rightarrow{1}$ can be taken to find the final solution. 
\subsection{Integral Transforms}
Using any of the known integral transforms, like Laplace and Fourier, has been the simplest way of finding solutions to infinite series with a parameter. If there is no suitable parameter, one can be added and at the end, set to $1$ or another value to make it disappear. The parameter must be located in the summand in such a functional position that the sum will still be uniformly convergent and can be Laplace transformed in either direction.  

We start with
\begin{equation}
s(\alpha)=\sum_{k=1}^{\infty}{g(k, \alpha)} \label{eqn22}
\end{equation}
and subject it to the inverse Laplace transform
\begin{equation}
\textsl{L}^{-1}_{\alpha}[\sum_{k=1}^{\infty}g(k, \alpha)],x=\sum_{k=1}^{\infty}G(k, x)=\textsl{L}^{-1}_{\alpha}[s(\alpha)],x \label{eqn23}
\end{equation}
To proceed, the elements of the sum must be transformable for this to be of any use. The transformed sum is attempted to be solved by any available methods. In a successful case, it will deliver the function
\begin{equation}
\textsl{L}^{-1}_{\alpha}[s(\alpha)],x=S(x) \label{eqn24}
\end{equation}
which can next be Laplace transformed to get the solution
\begin{equation}
\textsl{L}_{x}[S(x)],\alpha=s(\alpha) \label{eqn25}
\end{equation}
The validity of this approach is based on the validity of each operation in the process. Success of the method depends on the resulting parameter function and its transforms to make the sum solvable in equation (\ref{eqn23}). 

The Fourier transform can be used in a similar way and also the Fourier sine and Fourier cosine transforms can be used. The methods of integral transforms are somewhat ineffective and easily produce divergent series or give some formal solution which may turn out to be false. 

\subsection{Euler-Maclaurin Formula}
The Euler-Maclaurin formula \cite{Abramowitz1970}, \cite{Jeffrey2008}, \cite{Ivic1985} is the most important method for solving infinite series. It is applicable to sequences and generally used for finite difference problems. It gives a relationship between the sum and integrals as follows, assuming $a,b$ to be integers. Let the first $2n$ derivatives of $f(x)$ be continuous on an interval $[a,b]$. The interval is divided into equal parts $h=\frac{(a-b)}{n}$. Then for some $\theta , 0\leq{\theta}\leq{1}$ 
\begin{equation}
\sum_{k=0}^{m}{f(a+kh)}=\frac{1}{h}\int^{b}_{a}{f(x)dx}+\frac{1}{2}(f(a)+f(b))+\sum_{k=1}^{n-1}{B_{2k}h^{2k-1}\frac{[f^{(2k-1)}(b)-f^{(2k-1)}(a)]}{(2k)!}}+ \nonumber
\end{equation}
\begin{equation}
+\frac{h^{2n}}{(2n)!}{B_{2n}}\sum_{k=0}^{m-1}{f^{(2n)}(a+kh+\theta{h})} \label{eqn33}
\end{equation}
The $B_{i}$ are Bernoulli numbers. This formula is a bit tedious to use and the amount of work is dependent on the derivatives. A closed-form solution is rare and usually it gives only an approximation.

\subsection{Poisson Formula}
The Poisson formula usually does not give an immediate solution but is a transform allowing other procedures to be applied.
\begin{equation}
\sum_{k=-\infty}^{\infty}{F(k)}=\sum_{m=-\infty}^{\infty}\int^{\infty}_{-\infty}{e^{2\pi{imx}}F(x)dx} \label{eqn40}
\end{equation}

\subsection{Partial Summation}
In sum calculus partial summation is a useful tool \cite{Schaum1971}
\begin{equation}
\sum{y(k)\Delta{z(k)}}=y(k)z(k)-\sum{z(k+1)\Delta{y(k)}} \label{eqn50}
\end{equation}
where the $\Delta$ is the difference operator
\begin{equation}
\Delta{y(k)}=y(k+1)-y(k) \label{eqn60}
\end{equation}
This is analogous to integration by parts and may produce a solution.

\subsection{Power Series}
If we have a power series of the form (see \cite{Schaum1971})
\begin{equation}
\sum_{k=0}^{\infty}{a_k}{x^k}=\sum_{k=0}^{\infty}{u_k}{v_k}{x^k} \label{eqn70}
\end{equation}
where $u_k$ can be written  as a polynomial in k and $v_k$ is defined as
\begin{equation}
V(x)=\sum_{k=0}^{\infty}{v_k}{x^k} \label{eqn75}
\end{equation}
Then we have
\begin{equation}
\sum_{k=0}^{\infty}{a_k}{x^k}=V(x){u_0}+xV'(x)\frac{\Delta{u_0}}{1!}+{x^2}V''(x)\frac{\Delta^2{u_0}}{2!}+... \label{eqn78}
\end{equation}
Because $u_k$ is a polynomial, the series on the right will terminate. 

\subsection{A Trivial Summation}
A simple way for summation is the following \cite{Schaum1971}. 
We have a difference 
\begin{equation}
\Delta{u(k)}=u(k+1)-u(k)=f(k) \label{eqn80}
\end{equation}
We can subject it to the summation operator
\begin{equation}
\Delta^{-1}{\Delta{u_k}}=\Delta^{-1}{f(k)} \label{eqn90}
\end{equation}
which is equal to 
\begin{equation}
\sum_{k=1}^{N}{f(k)}=u(N+1)-u(1)  \label{eqn94} 
\end{equation}
This is the sum of the difference of the $u(k)$. 

\subsection{Composition}
In Section 2 we derive the method, its extension and explain the process of how to apply it. In Section 3 we illustrate the method with some example cases. Appendix A extends further the equivalences to better suit various cases. We have made the presentation very explicit avoiding mathematical rigor and detailed proofs. We use throughout lower case symbols ($h(s)$) for Laplace transformed functions and upper case ($H(t)$) for inverse transformed functions.

\section{Derivation of the Formulas}
\subsection{The Simple Form}
We have an infinite series
\begin{equation}
\sum_{k=1}^{\infty}{g(k)} \label{eqn100}
\end{equation}
to be solved with the assumption that it will converge. We assume the $g(k)$ to have an inverse Laplace transform $G(t)$
\begin{equation}
g(k)=\int^{\infty}_{0}{e^{-k{t}}{G(t)}dt}  \label{eqn110}
\end{equation}
The discrete argument $k\in{N^+}$ is a subset of the complex continuum. We can temporarily extend the range of the variable to $k\in{C}$, formally used in the Laplace transform. Thus we get
\begin{equation}
\sum_{k=1}^{\infty}{g(k)}=\sum_{k=1}^{\infty}\int^{\infty}_{0}{e^{-k{t}}{G(t)}dt} \label{eqn120}
\end{equation}
Since the sum is uniformly convergent we can interchange the order of summation and integration
\begin{equation}
\sum_{k=1}^{\infty}{g(k)}=\int^{\infty}_{0}{\sum_{k=1}^{\infty}e^{-k{t}}{G(t)}dt}=\int^{\infty}_{0}{dt{G(t)}\sum_{k=1}^{\infty}e^{-k{t}}} \label{eqn130}
\end{equation}
We know the sum to be equal to
\begin{equation}
\sum_{k=1}^{\infty}e^{-k{t}}=\frac{1}{e^{t}-1} \label{eqn150}
\end{equation}
and the series becomes
\begin{equation}
\sum_{k=1}^{\infty}{g(k)}=\int^{\infty}_{0}{\frac{dt{G(t)}}{e^{t}-1}} \label{eqn160}
\end{equation}
\newpage
In order for this to be true we insist on the following requirements:
\begin{itemize}
	\item 1. the series (\ref{eqn100}) converges
	\item 2. $g(k)$ has an inverse Laplace transform $G(t)$
	\item 3. the resulting integral converges
\end{itemize}
Failure to comply with these requirements will lead to either a diverging integral or a false solution. The series has been transformed to a rather simple integral. 

\subsection{The Extended Formulas}
We can parametrize equation (\ref{eqn160}) with $\alpha{\in{C}}$ and $\Re(\alpha)>0$ as follows
\begin{equation}
\sum_{k=1}^{\infty}{g(\alpha{k})}=\int^{\infty}_{0}{\frac{dt{G(t)}}{e^{\alpha{t}}-1}}=f(\alpha) \label{eqn200}
\end{equation}
Proof of this goes in the same way as with equation (\ref{eqn160}). Uniform convergence to $f(\alpha)$ is required of the series. Since $f(\alpha)$ is a function of the parameter, we may assume it to be the result of a Laplace transform. This is not necessarily so but at this point we can take it as an assumption. Subjecting this equation to an inverse Laplace transform gives
\begin{equation}
F(x)=\textsl{L}^{-1}_{\alpha}[f(\alpha)],x \label{eqn210}
\end{equation}
The left side of equation (\ref{eqn200}) is subjected to the same transform and we have
\begin{equation}
F(x)=\sum_{k=1}^{\infty}{\frac{G(\frac{x}{k})}{k}} \label{eqn220}
\end{equation}
Proof of this follows directly from the basic properties of the Laplace transform. These equations are equivalent, being either Laplace transforms or inverse transforms of each other. The new series in equation (\ref{eqn220}) must be uniformly convergent to $F(x)$, $x \in R, x \geq 0$.
\newpage
The schematic below shows the structure.
\begin{equation}
\int^{\infty}_{0}{\frac{dt{G(t)}}{e^{\alpha{t}}-1}}\  \  =\  \  \sum_{k=1}^{\infty}{g(\alpha{k})} \  \  =\  \  f(\alpha) \label{eqn230}
\end{equation}
\begin{equation}
       {\textsl{L}_{x},{\alpha}} \Uparrow \  \  \   \Downarrow{\textsl{L}^{-1}_{\alpha},{x}}   \nonumber
\end{equation}
\begin{equation}
\sum_{k=1}^{\infty}{\frac{G(\frac{x}{k})}{k}}\  \  =\  \  F(x) \label{eqn240}
\end{equation}
Any of these expressions can be converted to any of the others, having more possibilities for solving series than equation (\ref{eqn160}) alone does. The functions $F(x)$ and $f(\alpha)$ are solutions to the equations above. They can not be understood in the sense that any function could be expanded as a series using these equations. 

One should carefully note the distinction between the two ways of applying the Laplace transform, to the index function itself
\begin{equation}
g(k)=\int^{\infty}_{0}{e^{-k{t}}{G(t)}dt}=\textsl{L}_{t}[G(t)],k  \label{eqn245}
\end{equation}
and to the parameter-dependent functions
\begin{equation}
F(x)=\textsl{L}^{-1}_{\alpha}[f(\alpha)],x \label{eqn247}
\end{equation}
Using the functions $f(\alpha)$ and $F(x)$ poses the additional requirements of uniform convergence of pertinent series and that they are connected by the Laplace transform. 
Appendix A shows several modifications to equations (\ref{eqn230}) and (\ref{eqn240}). Each of them is proven in the same way.

\subsection{The Extended Method For Solving}
We may have a series of type A in a parametrized form 
\begin{equation}
\sum_{k=1}^{\infty}{g(\alpha{k})} \label{eqn310}
\end{equation}
or we might have a series similar to type B
\begin{equation}
\sum_{k=1}^{\infty}{\frac{G(\frac{x}{k})}{k}} \label{eqn312}
\end{equation}
\newpage
It is possible to proceed in different ways to solve a series:

	\begin{itemize}
  \item Type A,  equation (\ref{eqn310}), integral
	\begin{itemize}
	\item 1. starting from equation (\ref{eqn230})
	\item  2. determine the inverse Laplace transform $G(t)$
	\item  3. solve the integral in equation (\ref{eqn230}) and get $f(\alpha)$
	\item  4. the parameter $\alpha$ is set to $1$. If only equation (\ref{eqn230}) is applied, the parameter $\alpha$ is not necessary. The solution is $f(1)$
  \end{itemize}

  \item Type A,  equation (\ref{eqn310}), via type B series
	\begin{itemize}
	\item 1. starting from equation (\ref{eqn230})
  \item 2. calculate the inverse Laplace transform $G(t)$
  \item 3. solve the new series (\ref{eqn240}) getting $F(x)$ 
  \item 4. $F(x)$ is Laplace transformed to get $f(\alpha)$. The solution is $f(1)$	
  \end{itemize}

  \item Type B series equation (\ref{eqn312}), via type A series 
	\begin{itemize}
	\item 1. starting from equation (\ref{eqn240})
  \item 2. generate the Laplace transform $g(k)$ from $G(t)$
	\item 3. parametrize it to $g(\alpha{k})$ 
	\item 4. solve the sum in equation (\ref{eqn230})
	\item 5. inverse transform the resulting $f(\alpha)$ to get $F(x)$, the solution. 
	\item 6. Set $x$ to some final value depending on the fitting in equation (\ref{eqn240}).
  \end{itemize}

  \item Type B series equation (\ref{eqn312}), via integral and inverse transform
	\begin{itemize}
	\item 1. starting from equation (\ref{eqn240})
  \item 2. use the $G(t)$ to solve the integral in equation (\ref{eqn230}) to get $f(\alpha)$
	\item 3. inverse transform $f(\alpha)$ to get $F(x)$. 
	\item 4. Set $x$ to some final value depending on the fitting in equation (\ref{eqn240}).
  \end{itemize}
  \end{itemize}
Success of this method is dependent on the existence and eventual finding of Laplace transforms and inverse transforms necessary, on solving the new series or the resulting integral.

One is thus able to work in either space. This structure has an analogy with a solid state lattice having spatial vectors and its reciprocal lattice with wave vectors. 

\section{Example Cases}
In the following we display a few example applications of the method presented. We have systematically kept the $\alpha$ parameter in the equations though it might not be used. It would be necessary only if equation (\ref{eqn240}) would be required to be generated or the function $f(\alpha)$ itself would be of interest. We avoid proofs for simplicity.

\subsection{Riemann Zeta Function}
A simple example for applying equation (\ref{eqn230}) is the series of inverse power functions, with $\Re(z)>1$ 
\begin{equation}
\sum^{\infty}_{k=1}g({k})=\sum^{\infty}_{k=1}\frac{1}{{k}^z} \label{eqn500}
\end{equation}
and this is parametrized to
\begin{equation}
\sum^{\infty}_{k=1}g(\alpha{k})=\sum^{\infty}_{k=1}\frac{1}{(\alpha{k})^z}=f(\alpha) \label{eqn502}
\end{equation}
The inverse Laplace transform of $g(k)$ is
\begin{equation}
G(t)=\frac{t^{z-1}}{\Gamma(z)} \label{eqn510}
\end{equation}
Therefore we obtain by equation (\ref{eqn230}) 
\begin{equation}
\sum^{\infty}_{k=1}g(\alpha{k})=\int^{\infty}_{0}{\frac{dt{G(t)}}{e^{\alpha{t}}-1}}=\frac{1}{\Gamma(z)}\int^{\infty}_{0}{\frac{dt\cdot{t^{z-1}}}{e^{\alpha{t}}-1}} \label{eqn520}
\end{equation}
We recognize the integral representation of the Riemann $\zeta(z)$ and thus
\begin{equation}
\sum^{\infty}_{k=1}g(\alpha{k})=\frac{\zeta(z)}{\alpha^{z}} \label{eqn530}
\end{equation}
As $\alpha\rightarrow{1}$ we get finally
\begin{equation}
\sum^{\infty}_{k=1}g({k})=\zeta(z) \label{eqn540}
\end{equation}
closing the loop. To demonstrate equation (\ref{eqn240}) we turn our eyes to 
\begin{equation}
\sum_{k=1}^{\infty}{\frac{G(\frac{x}{k})}{k}}=\sum_{k=1}^{\infty}{\frac{x^{z-1}}{k^z{\Gamma(z)}}} \label{eqn550}
\end{equation}
\begin{equation}
=\frac{x^{z-1}\zeta(z)}{\Gamma(z)}=F(x) \nonumber
\end{equation}
Let's do Laplace on this and we will get
\begin{equation}
\textsl{L}_{x}[F(x)],{\alpha}=\frac{\zeta(z)}{\Gamma(z)}\int^{\infty}_{0}{dx\cdot{x^{z-1}}e^{-\alpha{x}}} \label{eqn560}
\end{equation}
\begin{equation}
=\frac{\zeta(z)}{\alpha^z}=f(\alpha) \nonumber
\end{equation}
getting the same result via another path.

\subsection{Simple Trigonometric Series}
The cosine series is a good test case since the outcome is simple. The parameter $\alpha$ is required to be $0<\alpha{<2\pi}$  
\begin{equation}
\sum^{\infty}_{k=1}cos(\alpha{k})=\sum^{\infty}_{k=1}g(\alpha{k})  \label{eqn600}
\end{equation}
The inverse Laplace transform of $g(k)$ is
\begin{equation}
G(t)=\frac{1}{2}[\delta(t+i)+\delta(t-i)] \label{eqn610}
\end{equation}
The integral in equation (\ref{eqn230}) gives us
\begin{equation}
\sum^{\infty}_{k=1}g(\alpha{k})=\int^{\infty}_{0}{\frac{dt{G(t)}}{e^{\alpha{t}}-1}}=\frac{1}{2}\int^{\infty}_{0}{\frac{dt\cdot{[\delta(t+i)+\delta(t-i)]}}{e^{\alpha{t}}-1}} \label{eqn620}
\end{equation}
Sorting out the integral produces
\begin{equation}
\sum^{\infty}_{k=1}g(\alpha{k})=-\frac{1}{2} \label{eqn630}
\end{equation}
irrespective of $\alpha$.

The integrated cosine series is a bit more complex. Therefore we use the new result in Appendix A. The parameter $\alpha$ is within $0<\alpha{<2\pi}$ 
\begin{equation}
\sum^{\infty}_{k=1}\frac{sin(\alpha{k})}{k}=\sum^{\infty}_{k=1}\frac{g(\alpha{k})}{k}  \label{eqn640}
\end{equation}
The inverse Laplace transform of $g(k)$ is
\begin{equation}
G(t)=\frac{1}{2i}[\delta(t+i)-\delta(t-i)] \label{eqn650}
\end{equation}
We use the integral in equation (\ref{eqn2530}) and obtain
\begin{equation}
\sum^{\infty}_{k=1}\frac{g(\alpha{k})}{k}=-\int^{\infty}_{0}{dt{G(t)ln(1-e^{-\alpha{t}})}} \nonumber
\end{equation}
\begin{equation}
=-\frac{1}{2i}\int^{\infty}_{0}{dt\cdot{ln(1-e^{-\alpha{t}})}\cdot{[\delta(t+i)-\delta(t-i)]}} \label{eqn660}
\end{equation}
The integral becomes
\begin{equation}
\sum^{\infty}_{k=1}\frac{sin(\alpha{k})}{k}=-\frac{1}{2i}ln(-e^{\alpha{i}})=\frac{\pm\pi}{2}-\frac{\alpha}{2} \label{eqn670}
\end{equation}
The positive signed solution is generally known. Numerical verification indicates that the negative sign is valid when $0>\alpha{>-2\pi}$. The equation should read
\begin{equation}
\sum^{\infty}_{k=1}\frac{sin(\alpha{k})}{k}=\frac{\pi\cdot{{signum(\alpha)}}}{2}-\frac{\alpha}{2} \label{eqn672}
\end{equation}
 As it happens, we keep the parameter.

\subsection{Converting a Fractional Series to a Power Series}
A more requiring case is the following convergent series, with $\left|a\right|<1$, $\Re(\beta)>1$ 
\begin{equation}
\sum^{\infty}_{k=1}g(k)=\sum^{\infty}_{k=1}\frac{1}{(a+k)^{\beta}} \label{eqn700}
\end{equation}
The inverse Laplace transform of $g(k)$ is
\begin{equation}
G(t)=\frac{t^{\beta-1}}{\Gamma(\beta)}e^{-a\cdot{t}} \label{eqn710}
\end{equation}
We have by equation (\ref{eqn230}) 
\begin{equation}
\sum^{\infty}_{k=1}g(\alpha{k})=\sum^{\infty}_{k=1}\frac{1}{(a+\alpha{k})^{\beta}}=\frac{1}{\Gamma(\beta)}\int^{\infty}_{0}{\frac{dt\cdot{e^{-at}t^{\beta-1}}}{e^{\alpha{t}}-1}} \label{eqn720}
\end{equation}
By expanding the exponential function to a power series and by using the integral representation of the Riemann $\zeta(s)$, $\Re(s)>1$ we get
\begin{equation}
\sum^{\infty}_{k=1}g(\alpha{k})=\sum^{\infty}_{n=0}\frac{(-1)^{n}a^{n}\Gamma(\beta+n)\zeta(\beta+n)}{\Gamma(n+1)\Gamma(\beta)\alpha^{\beta}} \label{eqn730}
\end{equation}
As $\alpha\rightarrow{1}$ we get 
\begin{equation}
\sum^{\infty}_{k=1}g(k)=\sum^{\infty}_{k=1}\frac{1}{(a+k)^{\beta}}=\sum^{\infty}_{n=0}\frac{(-1)^{n}a^{n}\Gamma(\beta+n)\zeta(\beta+n)}{\Gamma(n+1)\Gamma(\beta)} \label{eqn740}
\end{equation}
This is verified numerically with the ranges of validity mentioned above. The advantage of this result is that we have got the original less attractive function on the left expanded as a power series of $a$. On the other hand, this series is closely related to the Hurwitz zeta function.

\subsection{A Logarithmic Series Becomes a Zeta Function}
A more complicated case is the following alternating series, with $b,a \in C$, $\Re(b+ia+1) > 0$, $a\neq{0}$ and $\left|\Im{(a)}\right|<1$ 
\begin{equation}
\sum^{\infty}_{k=1}\frac{(-1)^{k+1}sin(a\cdot{ln(k)})}{k^{b+1}}=\sum^{\infty}_{k=1}{(-1)^{k+1}g(\alpha{k})} \label{eqn800}
\end{equation}
We plan to use equation (\ref{eqn980}) and need the inverse Laplace transform of
\begin{equation}
g(k)=\frac{sin(a\cdot{ln(k)})}{k^{b+1}} \label{eqn802}
\end{equation}
We can use the method for calculating inverse Laplace transforms for implicit logarithmic functions \cite{Stenlund2014} with $h(s)=\textsl{L}_{t}[H(t)],{s}$
\begin{equation}
\textsl{L}^{-1}_{s}[\frac{h(c\cdot{ln(s)})}{s^{b+1}}],{t}=\int^{\infty}_{0}{\frac{t^{cu+b}{H(u)}du}{\Gamma(cu+b+1)}} \label{eqn803}
\end{equation}
with
\begin{equation}
h(s)=sin(as) \label{eqn806}
\end{equation}
\begin{equation}
H(t)=\frac{1}{2i}[\delta{(t+ia)}-\delta{(t-ia)}] \label{eqn807}
\end{equation}
obtaining 
\begin{equation}
G(t)=\frac{1}{2i}[\frac{t^{b-ia}}{\Gamma(b-ia+1)}-\frac{t^{b+ia}}{\Gamma(b+ia+1)}] \label{eqn808}
\end{equation}
We use the basic integral definition of the Riemann zeta function for $\Re(s) > 0$
\begin{equation}
\zeta(s)=\frac{1}{1-2^{1-s}}\int^{\infty}_{0}{\frac{dt\cdot{t^{s-1}}}{e^t+1}} \label{eqn809}
\end{equation}
We have now all the information needed to apply equation (\ref{eqn980})  
\begin{equation}
\sum^{\infty}_{k=1}{(-1)^{k+1}g(\alpha{k})}=\frac{\alpha^{-b-1}}{2i}[\alpha^{ia}\zeta{(b-ia+1)}(1-2^{-b+ia})-\alpha^{-ia}\zeta{(b+ia+1)}(1-2^{-b-ia})] \label{eqn810}
\end{equation}
As $\alpha\rightarrow{1}$ we get 
\begin{equation}
\sum^{\infty}_{k=1}{(-1)^{k+1}g(k)}=\frac{1}{2i}[\zeta{(b-ia+1)}(1-2^{-b+ia})-\zeta{(b+ia+1)}(1-2^{-b-ia})] \nonumber
\end{equation}
\begin{equation}
=-\Im[{\zeta(b+ia+1)(1-2^{-b-ia})}] \label{eqn815}
\end{equation}
It is tempting to see how the companion function $cos(a\cdot{ln(k)})$ would behave. We process the equations through, starting from
\begin{equation}
\sum^{\infty}_{k=1}{(-1)^{k+1}m(k)}=\sum^{\infty}_{k=1}\frac{(-1)^{k+1}cos(a\cdot{ln(k)})}{k^{b+1}} \label{eqn820}
\end{equation}
The inverse Laplace is the following
\begin{equation}
M(t)=\frac{1}{2}[\frac{t^{b-ia}}{\Gamma(b-ia+1)}+\frac{t^{b+ia}}{\Gamma(b+ia+1)}] \label{eqn825}
\end{equation}
and we get
\begin{equation}
\sum^{\infty}_{k=1}{(-1)^{k+1}m(k)}=\frac{1}{2}[\zeta{(b-ia+1)}(1-2^{-b+ia})+\zeta{(b+ia+1)}(1-2^{-b-ia})] \nonumber
\end{equation}
\begin{equation}
=\Re[{\zeta(b+ia+1)(1-2^{-b-ia})}] \label{eqn830}
\end{equation}
If we now sum equations (\ref{eqn815}) and (\ref{eqn830}) properly, we get
\begin{equation}
\sum^{\infty}_{k=1}\frac{(-1)^{k+1}(cos(a\cdot{ln(k)})-i\cdot{sin(a\cdot{ln(k)})})}{k^{b+1}}=\sum^{\infty}_{k=1}\frac{(-1)^{k+1}}{k^{b+ia+1}} \nonumber
\end{equation}
\begin{equation}
=\zeta(b+ia+1)(1-2^{-b-ia}) \label{eqn840}
\end{equation}
The series expression on the right is equal to $\zeta(b+ia+1)(1-2^{-b-ia})$ when $\Re(b+ia+1) > 0$. Indeed, this was all about the $\zeta(s)$. 

\subsection{Partial Summation}
The equation (\ref{eqn160}) can be extended to a kind of partial summation if we apply it only to the function $f(k)$ in the following sum
\begin{equation}
\sum_{k=1}^{\infty}{g(k)}{f(k)}=\int^{\infty}_{0}{{dt{F(t)}}\sum_{k=1}^{\infty}{g(k)}{e^{-k{t}}}} \label{eqn850}
\end{equation}
Here the $F(t)$ is the inverse Laplace of $f(k)$
\begin{equation}
f(k)=\int^{\infty}_{0}{e^{-k{t}}{F(t)}dt}  \label{eqn860}
\end{equation}
The remaining summation may be more attractive for solving. The integral remains to be solved. 

\section{Discussion}
We have shown that the method presented offers a new way of generating solutions to infinite series. Various types of functions $g(k)$ of the index can be handled. The simple formula (\ref{eqn160}) is quite effective in many cases to offer a solution in the form of a rather simple integral. A parameter $\alpha$ can be added to the functions as a multiplier of the index in order to extend the method. The resulting equations (\ref{eqn230}) and (\ref{eqn240}) offer more options for getting solutions. They are equivalent to each other via a mediating Laplace transform. The method is explained in Chapter 2.3.

The conditions stated after equation (\ref{eqn160}) should be honored. It is essential for a new solution to be verified and the range of validity needs to be determined. Success of the method is dependent on finding the necessary Laplace forward or inverse transforms. Convergence or uniform convergence is required of the series. 


\newpage
\appendix
\section{Appendix. Various Particular Forms of Equivalence}
We can generate different equations modifying equation (\ref{eqn230}) and (\ref{eqn240}) to better suite practical applications. These forms are generally valid but they may be more useful than the basic formulas as their forms are already bent in the right direction. Proofs are similar to those before. We have added subscripts to the functions $f(\alpha)$ and $F(x)$ to remind that they differ between the cases. The formulas are given names referring to their expected crude behavior as a series. Notice the different index ranges in some cases. These equations can be used in the simple form of the upper formula with $\alpha=1$.

\subsection{Alternating Series}
\begin{equation}
\int^{\infty}_{0}{\frac{dt{G(t)}}{e^{\alpha{t}}+1}}\  \  =\  \  \sum_{k=1}^{\infty}{(-1)^{k+1}g(\alpha{k})} \  \  =\  \  f_0(\alpha) \label{eqn980}
\end{equation}
\begin{equation}
       {\textsl{L}_{x},{\alpha}} \Uparrow \  \  \   \Downarrow{\textsl{L}^{-1}_{\alpha},{x}}   \nonumber
\end{equation}
\begin{equation}
\sum_{k=1}^{\infty}{\frac{(-1)^{k+1}G(\frac{x}{k})}{k}}\  \  =\  \  F_0(x) \label{eqn985}
\end{equation}

\subsection{Shifted Series}
\begin{equation}
\int^{\infty}_{0}{\frac{dt{G(t)e^{-\beta{t}}}}{e^{\alpha{t}}-1}}\  \  =\  \  \sum_{k=1}^{\infty}{g(\alpha{k}+\beta)} \  \  =\  \  f_1(\alpha) \label{eqn990}
\end{equation}
\begin{equation}
       {\textsl{L}_{x},{\alpha}} \Uparrow \  \  \   \Downarrow{\textsl{L}^{-1}_{\alpha},{x}}   \nonumber
\end{equation}
\begin{equation}
\sum_{k=1}^{\infty}{\frac{e^{\frac{-\beta{x}}{k}}G(\frac{x}{k})}{k}}\  \  =\  \  F_1(x) \label{eqn1010}
\end{equation}

\subsection{Shifted Alternating Series}
\begin{equation}
\int^{\infty}_{0}{\frac{dt{G(t)e^{-\beta{t}}}}{e^{\alpha{t}}+1}}\  \  =\  \  \sum_{k=1}^{\infty}{(-1)^{k+1}g(\alpha{k}+\beta)} \  \  =\  \  f_2(\alpha) \label{eqn1050}
\end{equation}
\begin{equation}
       {\textsl{L}_{x},{\alpha}} \Uparrow \  \  \   \Downarrow{\textsl{L}^{-1}_{\alpha},{x}}   \nonumber
\end{equation}
\begin{equation}
\sum_{k=1}^{\infty}{\frac{e^{\frac{-\beta{x}}{k}}(-1)^{k+1}G(\frac{x}{k})}{k}}\  \  =\  \  F_2(x) \label{eqn1060}
\end{equation}

\subsection{Power Factor Series}
We have with $|{\gamma}|>1$
\begin{equation}
\int^{\infty}_{0}{\frac{dt{G(t)}}{\gamma{e^{\alpha{t}}}-1}}\  \  =\  \  \sum_{k=1}^{\infty}{\frac{g(\alpha{k})}{\gamma^k}} \  \  =\  \  f_3(\alpha) \label{eqn1100}
\end{equation}
\begin{equation}
       {\textsl{L}_{x},{\alpha}} \Uparrow \  \  \   \Downarrow{\textsl{L}^{-1}_{\alpha},{x}}   \nonumber
\end{equation}
\begin{equation}
\sum_{k=1}^{\infty}{\frac{G(\frac{x}{k})}{\gamma^k{k}}}\  \  =\  \  F_3(x) \label{eqn1110}
\end{equation}

\subsection{Power Factor Alternating Series}
With $|{\gamma}|>1$
\begin{equation}
\int^{\infty}_{0}{\frac{dt{G(t)}}{\gamma{e^{\alpha{t}}}+1}}\  \  =\  \  \sum_{k=1}^{\infty}{(-1)^{k+1}{\frac{g(\alpha{k})}{\gamma^k}}} \  \  =\  \  f_4(\alpha) \label{eqn1150}
\end{equation}
\begin{equation}
       {\textsl{L}_{x},{\alpha}} \Uparrow \  \  \   \Downarrow{\textsl{L}^{-1}_{\alpha},{x}}   \nonumber
\end{equation}
\begin{equation}
\sum_{k=1}^{\infty}{\frac{(-1)^{k+1}G(\frac{x}{k})}{\gamma^k{k}}}\  \  =\  \  F_4(x) \label{eqn1160}
\end{equation}

\subsection{Exponential Factor Series}
With $\Re(\beta)>0$
\begin{equation}
\int^{\infty}_{0}{\frac{dt{G(t)}}{e^{\alpha{t}+\beta}-1}}\  \  =\  \  \sum_{k=1}^{\infty}{e^{-\beta{k}}g(\alpha{k})} \  \  =\  \  f_5(\alpha) \label{eqn1180}
\end{equation}
\begin{equation}
       {\textsl{L}_{x},{\alpha}} \Uparrow \  \  \   \Downarrow{\textsl{L}^{-1}_{\alpha},{x}}   \nonumber
\end{equation}
\begin{equation}
\sum_{k=1}^{\infty}{\frac{e^{-\beta{k}}G(\frac{x}{k})}{k}}\  \  =\  \  F_5(x) \label{eqn1190}
\end{equation}

\subsection{Exponential Factor Alternating Series}
With $\Re(\beta)>0$
\begin{equation}
\int^{\infty}_{0}{\frac{dt{G(t)}}{e^{\alpha{t}+\beta}+1}}\  \  =\  \  \sum_{k=1}^{\infty}{(-1)^{k+1}{e^{-\beta{k}}g(\alpha{k})}} \  \  =\  \  f_6(\alpha) \label{eqn1210}
\end{equation}
\begin{equation}
       {\textsl{L}_{x},{\alpha}} \Uparrow \  \  \   \Downarrow{\textsl{L}^{-1}_{\alpha},{x}}   \nonumber
\end{equation}
\begin{equation}
\sum_{k=1}^{\infty}{\frac{(-1)^{k+1}e^{-\beta{k}}G(\frac{x}{k})}{k}}\  \  =\  \  F_6(x) \label{eqn1220}
\end{equation}

\subsection{Differentiated Series}
\begin{equation}
\int^{\infty}_{0}{\frac{dt{G(t)}}{(e^{\alpha{t}}-1)^2}}\  \  =\  \  \sum_{k=1}^{\infty}{k\cdot{g(\alpha{(k+1)})}} \  \  =\  \  f_7(\alpha) \label{eqn2330}
\end{equation}
\begin{equation}
       {\textsl{L}_{x},{\alpha}} \Uparrow \  \  \   \Downarrow{\textsl{L}^{-1}_{\alpha},{x}}   \nonumber
\end{equation}
\begin{equation}
\sum_{k=1}^{\infty}{\frac{{k}G(\frac{x}{k+1})}{k+1}}\  \  =\  \  F_7(x) \label{eqn2340}
\end{equation}

\subsection{Differentiated Alternating Series}
\begin{equation}
\int^{\infty}_{0}{\frac{dt{G(t)}}{(e^{\alpha{t}}+1)^2}}\  \  =\  \  \sum_{k=1}^{\infty}{(-1)^{k+1}k\cdot{g(\alpha{(k+1)})}} \  \  =\  \  f_8(\alpha) \label{eqn2430}
\end{equation}
\begin{equation}
       {\textsl{L}_{x},{\alpha}} \Uparrow \  \  \   \Downarrow{\textsl{L}^{-1}_{\alpha},{x}}   \nonumber
\end{equation}
\begin{equation}
\sum_{k=1}^{\infty}{\frac{(-1)^{k+1}G(\frac{x}{k+1})}{k+1}}\  \  =\  \  F_8(x) \label{eqn2440}
\end{equation}

\subsection{Integrated Series}
\begin{equation}
\int^{\infty}_{0}{dt{G(t)}{ln(1-e^{-\alpha{t}})}}\  \  =\  \  -\sum_{k=1}^{\infty}{\frac{g(\alpha{k})}{k}} \  \  =\  \  f_9(\alpha) \label{eqn2530}
\end{equation}
\begin{equation}
       {\textsl{L}_{x},{\alpha}} \Uparrow \  \  \   \Downarrow{\textsl{L}^{-1}_{\alpha},{x}}   \nonumber
\end{equation}
\begin{equation}
-\sum_{k=1}^{\infty}{\frac{G(\frac{x}{k})}{k^2}}\  \  =\  \  F_9(x) \label{eqn2540}
\end{equation}

\subsection{Integrated Alternating Series}
\begin{equation}
\int^{\infty}_{0}{dt{G(t)}{ln(1+e^{-\alpha{t}})}}\  \  =\  \  \sum_{k=1}^{\infty}{\frac{(-1)^{k+1}g(\alpha{k})}{k}} \  \  =\  \  f_{10}(\alpha) \label{eqn2630}
\end{equation}
\begin{equation}
       {\textsl{L}_{x},{\alpha}} \Uparrow \  \  \   \Downarrow{\textsl{L}^{-1}_{\alpha},{x}}   \nonumber
\end{equation}
\begin{equation}
\sum_{k=1}^{\infty}{\frac{(-1)^{k+1}G(\frac{x}{k})}{k^2}}\  \  =\  \  F_{10}(x) \label{eqn2640}
\end{equation}

\subsection{Added Constant Parameter Series}
We assume $\beta$ to be an extra constant parameter with $\left|\alpha\right|>\left|\beta\right|$
\begin{equation}
\int^{\infty}_{0}{\frac{dt{G(t)}}{e^{\alpha{t}}-e^{\beta{t}}}}\  \  =\  \  \sum_{k=1}^{\infty}{g(\alpha{k}-\beta{(k-1)})} \  \  =\  \  f_{11}(\alpha) \label{eqn2730}
\end{equation}
\begin{equation}
       {\textsl{L}_{x},{\alpha}} \Uparrow \  \  \   \Downarrow{\textsl{L}^{-1}_{\alpha},{x}}   \nonumber
\end{equation}
\begin{equation}
{e^{\beta{x}}}\sum_{k=1}^{\infty}{\frac{e^{\frac{-\beta{x}}{k}}{G(\frac{x}{k})}}{k}}\  \  =\  \  F_{11}(x) \label{eqn2740}
\end{equation}

\subsection{Alternating Added Constant Parameter Series}
We assume $\beta$ to be an extra constant parameter with $\left|\alpha\right|>\left|\beta\right|$
\begin{equation}
\int^{\infty}_{0}{\frac{dt{G(t)}}{e^{\alpha{t}}+e^{\beta{t}}}}\  \  =\  \  \sum_{k=1}^{\infty}{(-1)^{k+1}g(\alpha{k}-\beta{(k-1)})} \  \  =\  \  f_{12}(\alpha) \label{eqn2750}
\end{equation}
\begin{equation}
       {\textsl{L}_{x},{\alpha}} \Uparrow \  \  \   \Downarrow{\textsl{L}^{-1}_{\alpha},{x}}   \nonumber
\end{equation}
\begin{equation}
{e^{\beta{x}}}\sum_{k=1}^{\infty}{\frac{e^{\frac{-\beta{x}}{k}}{(-1)^{k+1}G(\frac{x}{k})}}{k}}\  \  =\  \  F_{12}(x) \label{eqn2760}
\end{equation}

\subsection{Hyperbolic Inverted Sine Series}
\begin{equation}
\int^{\infty}_{0}{\frac{dt{G(t)}}{e^{\alpha{t}}-e^{-\alpha{t}}}}\  \  =\  \  \sum_{k=1}^{\infty}{g(\alpha{(2k-1)})} \  \  =\  \  f_{13}(\alpha) \label{eqn2830}
\end{equation}
\begin{equation}
       {\textsl{L}_{x},{\alpha}} \Uparrow \  \  \   \Downarrow{\textsl{L}^{-1}_{\alpha},{x}}   \nonumber
\end{equation}
\begin{equation}
\sum_{k=1}^{\infty}{\frac{G(\frac{x}{2k-1})}{2k-1}}\  \  =\  \  F_{13}(x) \label{eqn2840}
\end{equation}

\subsection{Hyperbolic Inverted Cosine Series}
\begin{equation}
\int^{\infty}_{0}{\frac{dt{G(t)}}{e^{\alpha{t}}+e^{-\alpha{t}}}}\  \  =\  \  \sum_{k=1}^{\infty}{(-1)^{k+1}g(\alpha{(2k-1)})} \  \  =\  \  f_{14}(\alpha) \label{eqn2930}
\end{equation}
\begin{equation}
       {\textsl{L}_{x},{\alpha}} \Uparrow \  \  \   \Downarrow{\textsl{L}^{-1}_{\alpha},{x}}   \nonumber
\end{equation}
\begin{equation}
\sum_{k=1}^{\infty}{\frac{(-1)^{k+1}G(\frac{x}{2k-1})}{2k-1}}\  \  =\  \  F_{14}(x) \label{eqn2940}
\end{equation}

\subsection{Hyperbolic Inverted Sine Series With a Complex Argument}
\begin{equation}
\int^{\infty}_{0}{\frac{dt{G(t)}}{e^{{t(\beta+i\alpha)}}-e^{-{t(\beta+i\alpha)}}}}=\  \  \sum_{k=1}^{\infty}{g((\beta+i\alpha){(2k-1)})} \  \  =\  \  f_{15}(\alpha) \label{eqn3030}
\end{equation}
\begin{equation}
       {\textsl{L}_{x},{\alpha}} \Uparrow \  \  \   \Downarrow{\textsl{L}^{-1}_{\alpha},{x}}   \nonumber
\end{equation}
\begin{equation}
\sum_{k=1}^{\infty}{\frac{G(\frac{x}{i(2k-1)})e^{i{x}\beta}}{i(2k-1)}}\  \  =\  \  F_{15}(x)  \label{eqn3040}
\end{equation}

\subsection{Hyperbolic Inverted Cosine Series With a Complex Argument}
\begin{equation}
\int^{\infty}_{0}{\frac{dt{G(t)}}{e^{{t(\beta+i\alpha)}}+e^{-{t(\beta+i\alpha)}}}}=\  \  \sum_{k=1}^{\infty}{(-1)^{k+1}g((\beta+i\alpha){(2k-1)})} \  \  =\  \  f_{16}(\alpha) \label{eqn3130}
\end{equation}
\begin{equation}
       {\textsl{L}_{x},{\alpha}} \Uparrow \  \  \   \Downarrow{\textsl{L}^{-1}_{\alpha},{x}}   \nonumber
\end{equation}
\begin{equation}
\sum_{k=1}^{\infty}{\frac{(-1)^{k+1}G(\frac{x}{i(2k-1)})e^{i{x}\beta}}{i(2k-1)}}\  \  =\  \  F_{16}(x)  \label{eqn3140}
\end{equation}

\subsection{Hyperbolic Sine Series}
\begin{equation}
\int^{\infty}_{0}{dt{G(t)}sinh(e^{-\alpha{t}})}=\  \  \sum_{k=1}^{\infty}{\frac{g(\alpha{(2k-1)})}{(2k-1)!}} \  \  =\  \  f_{17}(\alpha) \label{eqn3230}
\end{equation}
\begin{equation}
       {\textsl{L}_{x},{\alpha}} \Uparrow \  \  \   \Downarrow{\textsl{L}^{-1}_{\alpha},{x}}   \nonumber
\end{equation}
\begin{equation}
\sum_{k=1}^{\infty}{\frac{G(\frac{x}{2k-1})}{(2k-1)(2k-1)!}}\  \  =\  \  F_{17}(x) \label{eqn3240}
\end{equation}

\subsection{Hyperbolic Cosine Series}
\begin{equation}
\int^{\infty}_{0}{dt{G(t)}cosh(e^{-\alpha{t}})}=\  \  g(0)+\sum_{n=1}^{\infty}{\frac{g(\alpha{2n})}{(2n)!}} \  \  =\  \  f_{18}(\alpha) \label{eqn3330}
\end{equation}
\begin{equation}
       {\textsl{L}_{x},{\alpha}} \Uparrow \  \  \   \Downarrow{\textsl{L}^{-1}_{\alpha},{x}}   \nonumber
\end{equation}
\begin{equation}
\delta(x){g(0)}+\sum_{n=1}^{\infty}{\frac{G(\frac{x}{2n})}{2n(2n)!}}\  \  =\  \  F_{18}(x) \label{eqn3340}
\end{equation}

\subsection{Square Root Series}
\begin{equation}
\int^{\infty}_{0}{\frac{dt{G(t)}}{\sqrt(1-{e^{-\alpha{t}}})}}\  \  =\  \  g(0)+\sum_{n=1}^{\infty}{\frac{1\cdot{3}\cdot{5}...(2n-1)g(\alpha{n})}{2\cdot{4}\cdot{6}...(2n)}} \  \  =\  \  f_{19}(\alpha) \label{eqn3430}
\end{equation}
\begin{equation}
       {\textsl{L}_{x},{\alpha}} \Uparrow \  \  \   \Downarrow{\textsl{L}^{-1}_{\alpha},{x}}   \nonumber
\end{equation}
\begin{equation}
\delta(x){g(0)}+\sum_{n=1}^{\infty}{\frac{1\cdot{3}\cdot{5}...(2n-1)G(\frac{x}{n})}{2\cdot{4}\cdot{6}...(2n)n}}\  \  =\  \  F_{19}(x) \label{eqn3440}
\end{equation}

\subsection{Alternating Square Root Series}
\begin{equation}
\int^{\infty}_{0}{\frac{dt{G(t)}}{\sqrt(1+{e^{-\alpha{t}}})}}\  \  =\  \  g(0)+\sum_{n=1}^{\infty}{(-1)^{n}\frac{1\cdot{3}\cdot{5}...(2n-1)g(\alpha{n})}{2\cdot{4}\cdot{6}...(2n)}} \  \  =\  \  f_{20}(\alpha) \label{eqn3530}
\end{equation}
\begin{equation}
       {\textsl{L}_{x},{\alpha}} \Uparrow \  \  \   \Downarrow{\textsl{L}^{-1}_{\alpha},{x}}   \nonumber
\end{equation}
\begin{equation}
\delta(x){g(0)}+\sum_{n=1}^{\infty}{(-1)^{n}\frac{1\cdot{3}\cdot{5}...(2n-1)G(\frac{x}{n})}{2\cdot{4}\cdot{6}...(2n)n}}\  \  =\  \  F_{20}(x) \label{eqn3540}
\end{equation}

\subsection{Exponential Series}
\begin{equation}
\int^{\infty}_{0}{dt{G(t)}e^{e^{-\alpha{t}}}}\  \  =\  \  g(0)+\sum_{k=1}^{\infty}{\frac{g(\alpha{k})}{k!}} \  \  =\  \  f_{21}(\alpha) \label{eqn3550}
\end{equation}
\begin{equation}
       {\textsl{L}_{x},{\alpha}} \Uparrow \  \  \   \Downarrow{\textsl{L}^{-1}_{\alpha},{x}}   \nonumber
\end{equation}
\begin{equation}
\delta{(x)}g(0)+\sum_{k=1}^{\infty}{\frac{G(\frac{x}{k})}{k\cdot{k!}}}\  \  =\  \  F_{21}(x) \label{eqn3560}
\end{equation}

\subsection{Negative Exponential Series}
\begin{equation}
\int^{\infty}_{0}{dt{G(t)}e^{-e^{-\alpha{t}}}}\  \  =\  \  g(0)+\sum_{k=1}^{\infty}{\frac{(-1)^{k}g(\alpha{k})}{k!}} \  \  =\  \  f_{22}(\alpha) \label{eqn3570}
\end{equation}
\begin{equation}
       {\textsl{L}_{x},{\alpha}} \Uparrow \  \  \   \Downarrow{\textsl{L}^{-1}_{\alpha},{x}}   \nonumber
\end{equation}
\begin{equation}
\delta{(x)}g(0)+\sum_{k=1}^{\infty}{\frac{(-1)^{k}G(\frac{x}{k})}{k\cdot{k!}}}\  \  =\  \  F_{22}(x) \label{eqn3580}
\end{equation}

\end{document}